\documentclass{amsart}

\usepackage{amsfonts,amsmath,amssymb,amsthm,amscd,amsxtra}
\usepackage{enumerate,verbatim}
\usepackage{mathptmx}
\usepackage[all,2cell,ps]{xy}
\usepackage[pagebackref]{hyperref}

\theoremstyle{plain} %default

\newtheorem{thm}{Theorem}[section]

\newtheorem{prop}[thm]{Proposition}
\newtheorem{cor}[thm]{Corollary}

\theoremstyle{definition}
\newtheorem{lem}[thm]{Lemma}
\newtheorem{dfn}[thm]{Definition}

\newtheorem{rmk}[thm]{Remark}

\numberwithin{equation}{section}
\newcommand{\fa}{\mathfrak{a}}
\newcommand{\fb}{\mathfrak{b}}

\newcommand{\fm}{\mathfrak{m}}

\newcommand{\fc}{\mathfrak{c}}

\newcommand{\ZZ}{\mathbb{Z}}

\DeclareMathOperator{\E}{E}
\DeclareMathOperator{\ann}{ann} 
 
\DeclareMathOperator{\CI}{\textnormal{CI-dim}}

\DeclareMathOperator{\cx}{cx} \DeclareMathOperator{\Tr}{Tr}
\DeclareMathOperator{\gr}{\textnormal{grade}}
\DeclareMathOperator{\hte}{\textnormal{ht}}
\DeclareMathOperator{\G-dim}{\textnormal{G-dim}}

\DeclareMathOperator{\coker}{coker} 
\DeclareMathOperator{\depth}{depth} 
 \DeclareMathOperator{\Ext}{Ext}

\DeclareMathOperator{\g}{G}
\DeclareMathOperator{\hh}{H}
\DeclareMathOperator{\Hom}{Hom} 
\DeclareMathOperator{\id}{id}

\DeclareMathOperator{\pd}{pd} 
 
 \DeclareMathOperator{\rank}{rank}

\DeclareMathOperator{\Tor}{Tor} 
 
\DeclareMathOperator{\cm}{CM}

\def\urltilda{\kern -.15em\lower .7ex\hbox{\~{}}\kern .04em}
\def\urldot{\kern -.10em.\kern -.10em}\def\urlhttp{http\kern -.10em\lower -.1ex
\hbox{:}\kern -.12em\lower 0ex\hbox{/}\kern -.18em\lower
0ex\hbox{/}}

\begin{document}

\title[Symmetry in vanishing of Tate cohomology over Gorenstein rings]{Symmetry in vanishing of Tate cohomology over Gorenstein rings}

\author[A. Sadeghi]{Arash Sadeghi}
\address{School of Mathematics, Institute for Research in Fundamental Sciences (IPM), P.O. Box: 19395-5746, Tehran, Iran}
 \email{sadeghiarash61@gmail.com}


\keywords{complete intersection dimension, complexity, vanishing of Tate cohomology, linkage of modules.\\
This research was supported by a grant from IPM}

%\date{\today}

\begin{abstract}
We investigate symmetry in the vanishing of Tate cohomology for finitely generated modules over local Gorenstein
rings. For finitely generated $R$--modules $M$ and $N$ over Gorenstein local ring $R$, it is shown that $\widehat{\Ext}^i_R(M,N)=0$ for all $i\in\ZZ$ if and only if $\widehat{\Ext}^i_R(N,M)=0$ for all $i\in\ZZ$.
\end{abstract}

\maketitle{}

%\setcounter{tocdepth}{1}
%\tableofcontents
%%%%%%%%%%%%%%%%%%%%%%%%%%%%%%%%%%%%%%%%%%%%%%%%%%%%%%%%%%%%%%%%%%%%%%%%%%%%%%%%%%%%%%%%%%%%%%%%%%%%%%%%%%%%
\section{Introduction}
The symmetry in the vanishing of Ext for finitely generated modules
was studied by Avramov and Buchweitz \cite{AvBu}. For finitely generated $R$-modules $M$ and $N$ over
complete intersection ring $R$, they proved that the following are equivalent:
\begin{enumerate}[(i)]
\item{$\Ext^i_R(N,M)=0$ for all $i\gg0$.}
\item{$\Ext^i_R(M,N)=0$ for all $i\gg0$.}
\item{$\Tor_i^R(M,N)=0$ for all $i\gg0$.}
\end{enumerate}
In \cite{HJ}, Huneke and Jorgensen proved that the symmetry in vanishing of Ext holds for a wider class of Gorenstein rings. More precisely, they introduced a new class of rings between complete
intersection and Gorenstein rings, which is called AB ring. An AB
ring $R$ is a local Gorenstein ring defined by the property that
there is a constant $C$, depending only on the ring, such that if
$\Ext^i_R(M,N)=0$ for all $i\gg0$, then $\Ext^i_R(M,N)=0$ for all
$i> C$. For finitely
generated modules $M$ and $N$ over an AB ring $R$, they proved that $\Ext^i_R(M,N)=0$
for all $i\gg0$ if and only if $\Ext^i_R(N,M)=0$ for all $i\gg0$ $($see \cite[Theorem 4.1]{HJ}). The class of AB rings is strictly contained between complete intersections and Gorenstein rings (see \cite[Theorem 3.6]{HJ} and \cite[Theorem]{JS}). In \cite[Theorem 1.2]{JS1}, Jorgensen and \c{S}ega proved that the symmetry in vanishing of Ext does not hold for all Gorenstein rings.

In this paper, we investigate symmetry in the vanishing of Tate cohomology for finitely generated modules over Gorenstein local rings.
As a consequence, it is shown that there is a different kind of symmetry over Gorenstein rings. More precisely, we obtain the following result (see Theorem \ref{t4}).
\begin{thm}\label{tet}
	Let $R$ be a Gorenstein local ring of dimension $d$ and let $M$, $N$ be $R$--modules. Assume that $n$ is an integer.
	Then the following statements hold:
	\begin{enumerate}[(i)]
		\item{If $\widehat{\Ext}^i_R(M,N)=0$ for all $i\geq n$, then $\widehat{\Ext}^i_R(N,M)=0$ for all $i<-n$.}
		\item{If $\widehat{\Ext}^i_R(M,N)=0$ for all $i<n$, then $\widehat{\Ext}^i_R(N,M)=0$ for all $i\geq d-n$.}
	\end{enumerate}
\end{thm}
In \cite[Proposition 3.2]{CJ}, for finitely
generated modules $M$ and $N$ over an AB ring $R$, it is shown that $\Ext^i_R(M,N)=0$ for $i\gg0$ if and only if $\widehat{\Ext}^i_R(M,N)=0$ for all $i\in\ZZ$. Therefore, our main result recovers Theorem of 
Huneke and Jorgensen.

The organization of the paper is as follows. In Section 2, we collect preliminary notions, definitions and some known results which will be used in this paper.

In Section 3, we study the vanishing of Tate cohomology for finitely generated modules over Gorenstein local rings. As a consequence, we prove Theorem \ref{tet} in this Section.
We also generalize the Auslander-Reiten duality for modules of finite Gorenstein dimension (see Theorem \ref{t1}).
\begin{thm}
Let $(R,\fm,k)$ be a Cohen-Macaulay local ring of dimension $d$ with canonical module $\omega_R$. Assume that $M$, $N$ are $R$--modules of finite Gorenstein dimension and that $\widehat{\Ext}^i_R(M,N)$ has finite length for all $i\in\ZZ$. Then $$\widehat{\Ext}^i_R(M,N)\cong\Hom_R(\widehat{\Ext}^{d-i-1}_R(N,M\otimes_R\omega_R),\E_R(k))$$
for all $i\in\ZZ$, where $\E_R(k)$ is the injective envelope of the residue field of $R$.
\end{thm}
In \cite[Thoerem 4.1]{J}, J${\o}$rgensen proved that the symmetry in the vanishing of Ext holds for modules of finite complete intersection dimension over local Gorenstein rings. The concept of modules with reducible complexity was introduced by Bergh. Every module of finite complete intersection dimension has reducible complexity. In \cite[Theorem 3.5]{B2}, Bergh generalized the J${\o}$rgensen's result for modules with reducible complexity.
In Section 4, as a consequence of our main theorem, it is shown that the symmetry in the vanishing of Ext holds for modules over local Gorenstein rings, provided at least one of them has reducible complexity (see Corollary \ref{c7}).

In Section 5, we study the vanishing of Tate homology for $\g$-perfect $R$--modules.
Let $M$ be a perfect $R$--module and $N$ an $R$--module. The following equality is well-known:
$$\sup\{i\geq0\mid\Tor_i^R(M,N)\neq0\}+\inf\{i\geq0\mid\Ext^i_R(M,N)\neq0\}=\pd_R(M).$$
In Theorem \ref{t3}, we generalize the above equality for $G$-perfect modules.

Recall, an ideal $I$ in $R$ is said to be a Gorenstein ideal if $R/I$ is a Gorenstein local ring. As another application of our main result, we have the following (see Corollary \ref{c2}).
\begin{cor}\label{cc}
Let $R$ be a Gorenstein ring and let $\fa$, $\fb$ be Gorenstein ideals of $R$. The following conditions are equivalent:
\begin{enumerate}[(i)]
\item{$\Ext^i_R(R/\fa,R/\fb)=0$ for $i\gg0$;}
\item{$\Ext^i_R(R/\fb,R/\fa)=0$ for $i\gg0$;}
\item{$\Tor_i^R(R/\fa,R/\fb)=0$ for $i\gg0$;}
\end{enumerate}
\end{cor}
Finally, in the last section, we study the vanishing of Tate (co)homology for linked modules. In Proposition \ref{t44} and Theorem \ref{t6}, we generalize some results of Puthenpurakal \cite[Theorem 1.2 and Corollary 1.3]{Pu}.
As a consequence, we obtain the following surprising result:
\begin{cor}
Let $R$ be a complete Kleinian singularity and let $M$, $N$ be $R$--modules.
Then $\widehat{\Ext}^i_R(M,N)\cong\widehat{\Ext}^{i}_R(N,M)$  for all  $i\in\ZZ$. In particular, if $M$ and $N$ are maximal Cohen-Macaulay, then
$\Ext^i_R(M,N)\cong\Ext^{i}_R(N,M)$  for all  $i>0$.
\end{cor}
We also investigate the relationship between Betti and Bass numbers of linked module. As a consequence, it is shown that the stable Betti and Bass numbers are preserved under evenly linkage. More precisely, let $R$ be a Gorenstein local ring and let $\fa$ and $\fb$ be perfect Gorenstein ideals of $R$. Assume that $M$, $N$ and $L$ are Cohen-Macaulay $R$--modules such that $M\underset{\fa}{\sim}N\underset{\fb}{\sim}L$. In Theorem \ref{t2}, it is shown that $\widehat{\beta}_i^R(M)=\widehat{\beta}_i^R(L)$ and $\widehat{\mu}^i_R(M)=\widehat{\mu}^i_R(L)$ for all $i\in\ZZ$.
%%%%%%%%%%%%%%%%%%%%%%%%%%%%%%%%%%%%%%%%%%%%%%%%%%%%%%%%%%%%%%%%%%%%%%%%%%%%%%%%%%%%%%%%%%%%%%%%%%%%%%%%%%%%%
\section{Preliminaries}
Throughout the paper, $R$ is a commutative Noetherian local ring and all modules are finite (i.e. finitely generated)
$R$--modules. Let $$\cdots\rightarrow F_{n+1} \rightarrow F_{n}\rightarrow F_{n-1}\rightarrow\cdots\rightarrow F_0\rightarrow M\rightarrow0$$ be the minimal free resolution of $M$. Recall that the $n$th syzygy of an $R$--module $M$ is the cokernel of the $F_{n+1}\rightarrow F_{n}$ and denoted by $\Omega^n_R(M)$, and it is unique up to isomorphism. The $n$th Betti number, denoted $\beta_n^R(M)$, is the rank of the free $R$--module $F_n$.
The complexity of $M$ is defined as follows;
$$\cx_R(M)=\inf\{i\in \mathbb{N} \cup 0\mid \exists \gamma\in \mathbb{R}  \text{ such that } \beta_n^R(M)\leq\gamma n^{i-1}\text{ for } n\gg0\}.$$
Note that $\cx_R(M)=\cx_R(\Omega^i_R(M))$ for every $i\geq0$. It follows from the definition that
$\cx_R(M)=0$ if and only if $\pd_R(M)<\infty$.
The complete intersection dimension was introduced by Avramov, Gasharov and Peeva \cite{AGP}. A module of finite complete intersection dimension behaves homologically like a module over a complete intersection.
Recall that a quasi-deformation of $R$ is a diagram $R\rightarrow A\twoheadleftarrow Q$ of local homomorphisms, in which
$R\rightarrow A$ is faithfully flat, and $A\twoheadleftarrow Q$ is surjective with kernel generated by a regular sequence.
The module $M$ has finite complete intersection dimension if there exists such a quasi-deformation for which $\pd_Q(M\otimes_RA)$ is finite.
The complete intersection dimension of $M$, denoted $\CI_R(M)$, is defined as follows;
$$\CI_R(M)=\inf\{\pd_Q(M\otimes_RA)-\pd_Q(A)\mid R\rightarrow A\twoheadleftarrow Q \text{ is a quasi-deformation }\}.$$
Note that every module of finite complete intersection dimension has finite complexity \cite[Theorem 5.3]{AGP}.

The concept of modules with reducible complexity was introduced by Bergh \cite{B2}.\\
Let $M$ and $N$ be $R$--modules and consider a homogeneous element $\eta$ in the graded $R$--module
$\Ext^*_R(M,N)=\bigoplus^{\infty}_{i=0}\Ext^i_R(M,N)$. Choose a map $f_{\eta}:\Omega^{|\eta|}_R(M)\rightarrow N$ representing $\eta$, and denote by $K_{\eta}$ the pushout of this map and the inclusion $\Omega^{|\eta|}_R(M)\hookrightarrow F_{|\eta|-1}$.
Therefore we obtain a commutative diagram
$$\begin{CD}
&&&&&&&&\\
  \ \ &&&&  0@>>>\Omega^{|\eta|}_R (M) @>>> F_{|\eta|-1}@>>>\Omega^{|\eta|-1}_R (M) @>>>0&  \\
                                &&&&&& @VV{f_{\eta}}V @VV V @VV{\parallel} V\\
  \ \  &&&& 0@>>> N @>>> K_{\eta} @>>>\Omega^{|\eta|-1}_R (M)@>>>0.&\\
\end{CD}$$\\
with exact rows. Note that the module $K_{\eta}$ is independent, up to isomorphism, of the map $f_{\eta}$ chosen to represent ${\eta}$.
\begin{dfn}\cite{B2}\emph{
The full subcategory of $R$-modules consisting of the modules having
reducible complexity is defined inductively as follows:}
\begin{itemize}
        \item[(i)]\emph{ Every $R$-module of finite projective dimension has reducible complexity.}
         \item[(ii)]\emph{ An $R$-module $M$ of finite positive complexity has reducible complexity if
                    there exists a homogeneous element $\eta\in\Ext^{*}_R(M,M)$, of positive degree,
                         such that $\cx_R(K_{\eta}) < \cx_R(M)$, $\depth_R(M)=\depth_R(K_{\eta})$ and $K_{\eta}$ has reducible complexity.}
\end{itemize}
\end{dfn}
By \cite[Proposition 2.2(i)]{B2}, every module of finite complete intersection dimension has reducible complexity. On the other hand, there are modules having reducible complexity but whose complete intersection dimension is infinite (see for example, \cite[Corollary 4.7]{BJ}).

The notion of the Gorenstein(or G-) dimension was introduced by Auslander \cite{A1}, and developed by Auslander and Bridger in \cite{AB}.
\begin{dfn}
An $R$--module $M$ is said to be of $G$-dimension zero (or totally reflexive) whenever
\begin{itemize}
            \item[(i)]{\emph{the biduality map $M\rightarrow M^{**}$ is an isomorphism;}}
            \item[(ii)]{\emph{$\Ext^i_R(M,R)=0$ for all $i>0$;}}
            \item[(iii)]{\emph{$\Ext^i_R(M^*,R)=0$ for all $i>0$.}}
\end{itemize}
\end{dfn}
The Gorenstein dimension of $M$, denoted $\G-dim_R(M)$, is defined to be the infimum of all
nonnegative integers $n$, such that there exists an exact sequence
$$0\rightarrow G_n\rightarrow\cdots\rightarrow G_0\rightarrow  M \rightarrow 0$$
in which all the $G_i$ have $\g$-dimension zero. Note that $\G-dim_R(M)$ is bounded above by the complete intersection dimension, $\CI_R(M)$, of $M$ and if $\CI_R(M)<\infty$ then the equality holds (see \cite[Theorem 1.4]{AGP}).

For a finite presentation $P_1\overset{f}{\rightarrow}P_0\rightarrow
M\rightarrow 0$ of an $R$--module $M$, its transpose, $\Tr M$, is
defined as $\coker f^*$, where $(-)^* := \Hom_R(-,R)$, which
satisfies in the exact sequence
\begin{equation}\label{1.1}
0\rightarrow M^*\rightarrow P_0^*\rightarrow P_1^*\rightarrow \Tr
M\rightarrow 0.
\end{equation}
Moreover, $\Tr M$ is unique up to projective equivalence. Thus all
minimal projective presentations of $M$ represent isomorphic
transposes of $M$.

In the following, we summarize some basic facts about Gorenstein dimension (see \cite{AB} for more details).
\begin{thm}\label{G}
For an $R$--module $M$, the following statements hold.
\begin{enumerate}[(i)]
       \item{$\G-dim_R(M)=0$ if and only if $\G-dim_R(\Tr M)=0$;}
       \item{If $\G-dim_R(M)<\infty$ then $\G-dim_R(M)=\sup\{i\geq0\mid\Ext^i_R(M,R)\neq0\}$;}
        \item{If $\G-dim_R(M)<\infty$, then $\G-dim_R(M)=\depth R-\depth_R(M)$;}
        \item{$R$ is Gorenstein if and only if $\G-dim_R(M)<\infty$ for all finitely generated $R$--module $M$.}
\end{enumerate}
\end{thm}
Two modules $M$ and $N$ are called \emph{stably isomorphic} and write $M\approx N$ if $M\oplus P\cong N\oplus Q$
for some projective modules $P$ and $Q$.
%%%%%%%%%%%%%%%%%%%%%%%%%%%%%%%%%%%%%%%%%%%%%%%%%%%%%%%%%%%%%%%%%%%%%%%%%
\begin{thm}\cite[Theorem 2.8]{AB}\label{a1}
Let $M$ be an $R$--module and $n\geq0$ an integer. Then there are exact sequences of functors:
$$\Ext^1_R(\Tr\Omega^nM,-)\hookrightarrow\Tor_n^R(M,-)\rightarrow\Hom_R(\Ext^n_R(M,R),-)
\rightarrow\Ext^2_R(\Tr\Omega^{n}M,-),$$
$$\Tor_2^R(\Tr\Omega^{n}M,-)\rightarrow(\Ext^n_R(M,R)\otimes_R-)\rightarrow\Ext^n_R(M,-)\twoheadrightarrow
\Tor_1^R(\Tr\Omega^{n}M,-).$$
\end{thm}
Tate cohomology for modules of finite Gorenstein dimension was studied by Avramov and Martsinkovsky in \cite{AM}. Let $\mathcal{P}=\mathcal{P}(R)$ denote the full subcategory of $\mathcal{C}$(R) whose objects are the finite projective $R$--modules, where $\mathcal{C}(R)$ denotes the category of all $R$--module and $R$--homomorphisms.

A complex $\mathbf{T}$ is totally acyclic if $T_n\in\mathcal{P}$ and $\hh_n(\mathbf{T})=0=\hh_n(\Hom_R(\mathbf{T},R))$ for all $n\in\ZZ$.
A \emph{complete resolution} of an $R$-module $M$ is a diagram
$$ \mathbf{T} \overset{\vartheta}{\longrightarrow} \mathbf{P} \overset{\pi}{\longrightarrow}M,$$
where $\pi$ is a $\mathcal{P}$-resolution, $T$ is a totally
acyclic complex, $\vartheta$ is a morphism, and $\vartheta_i$ is an isomorphism
for $i \gg 0$. An $R$-module has finite $\g$-dimension if and only
if it has a complete resolution \cite[Theorem 3.1]{AM}.

Let $M$ be an $R$-module with a complete resolution $\mathbf{T}\rightarrow \mathbf{P}\rightarrow M$. For an $R$-module $N$, \emph{Tate cohomology} of $M$ and $N$ is defined as
$$ \widehat{\Ext}^i_R(M,N) = \hh^i(\Hom_R(\mathbf{T},N))  \text{ for } i\in\ZZ.$$
Also, for an $R$-module $N$, \emph{Tate homology} of $M$ and $N$ is defined as
$$\widehat{\Tor}_i^R(M,N) = \hh_i(\mathbf{T}\otimes_RN)  \text{ for } i\in\ZZ.$$
By construction, there are isomorphisms
\begin{equation}\label{t}
\widehat{\Tor}_i^R(M,N)\cong\Tor_i^R(M,N) \text{ and } \widehat{\Ext}^i_R(M,N)\cong\Ext^i_R(M,N),
\end{equation}
for all $i>\G-dim_R(M)$.

In the following we summarize some basic properties about Tate cohomology which will be used throughout the paper; see \cite{AM}, \cite{AvBu} and \cite{CJ1} for more details.
\begin{thm}\label{tate}Let $M$ be an $R$--module of finite Gorenstein dimension. The followings hold:
	 	\begin{enumerate}[(i)]
		\item{If $\pd_R(M)<\infty$ or $\pd_R(N)<\infty$, then $\widehat{\Ext}^i_R(M,N)=0$ for all $i\in\ZZ$.}
		\item{If $M$ is totally reflexive then $\widehat{\Tor}_i^R(M,N)\cong\widehat{\Ext}^{-i-1}_R(M^*,N)$ for all $i\in\mathbb{Z}$.}
		\item{If $N$ has finite injective dimension, then $\widehat{\Ext}^i_R(M,N)=0$ for all $i\in\ZZ$.}
		\item{For every exact sequence $0\rightarrow X\rightarrow Y\rightarrow Z\rightarrow0$ of $R$--modules, there is
			a doubly infinite long exact sequence $$\cdots\rightarrow\widehat{\Ext}^i_R(M,X)\rightarrow\widehat{\Ext}^i_R(M,Y)\rightarrow\widehat{\Ext}^i_R(M,Z)\rightarrow\cdots.$$}
			\item{For every exact sequence $0\rightarrow L\rightarrow M\rightarrow N\rightarrow0$
			of $R$-modules of finite $G$-dimension and each $R$--module $X$ there is
				a doubly infinite long exact sequence $$\cdots\rightarrow\widehat{\Ext}^i_R(N,X)\rightarrow\widehat{\Ext}^i_R(M,X)\rightarrow\widehat{\Ext}^i_R(L,X)\rightarrow\cdots.$$}
	\end{enumerate}
\end{thm}
\begin{rmk}\label{rem}
Let $M$ and $N$ be $R$--modules of finite Gorenstein dimension. Consider the following exact sequences:
\begin{equation}\tag{\ref{rem}.1}
	0\rightarrow X\rightarrow G\rightarrow M\rightarrow0 \text{ and }	
	0\rightarrow Y\rightarrow H\rightarrow N\rightarrow0,
	\end{equation}
where $\G-dim_R(G)=0=\G-dim_R(H)$ and $X$, $Y$ have finite projective dimension. By Theorem \ref{tate}(i), $\widehat{\Ext}^i_R(-,Y)=0=\widehat{\Ext}^i_R(-,X)$ for all $i\in\ZZ$. Hence, by Thoerem \ref{tate}(iv), (v), the exact sequence (\ref{rem}.1) induces the following isomorphisms
$$\widehat{\Ext}^i_R(M,N)\cong\widehat{\Ext}^i_R(G,H) \text{ and } \widehat{\Ext}^i_R(N,M)\cong\widehat{\Ext}^i_R(H,G),$$ for all $i\in\ZZ$. 
Also, by \cite[Lemma 2.7]{CJ1}, $\widehat{\Tor}_i^R(-,Y)=0=\widehat{\Tor}_i^R(-,X)$ for all $i\in\ZZ$.  Hence, it follows from \cite[Proposition 2.8, 2.9]{CJ1} that
$\widehat{\Tor}_i^R(M,N)\cong\widehat{\Tor}_i^R(G,H)$ for all $i\in\ZZ$. 
\end{rmk}
%%%%%%%%%%%%%%%%%%%%%%%%%%%%%%%%%%%%%%%%%%%%%%%%%%%%%%%%%%%%%%%%%%%%%%%%%%%%%%%%%%%%%%%%%%%%%%%%%%%%%%%%%

%%%%%%%%%%%%%%%%%%%%%%%%%%%%%%%%%%%%%%%%%%%%%%%%%%%%%%%%%%%%%%%%%%%%%%%%%%%%%%%%%%%%%
\section{Vanishing of Tate Cohomology}
In this section, we investigate the symmetry in the vanishing of Tate cohomology for modules over Gorenstein local rings.
For the proof of our main result we need
the following preliminary lemmas.
\begin{lem}\label{l2}
Let $(R,\fm,k)$ be a local ring and let $M$ and $N$ be $R$--modules. If $M$ has finite Gorenstein dimension, then
$$\widehat{\Ext}^i_R(M,N^{\vee})\cong
\widehat{\Tor}_i^R(M,N)^{\vee} \text{ for all } i\in\ZZ,$$
where $(-)^\vee:=\Hom_R(-,\E_R(k)).$
\end{lem}
\begin{proof}
Let $\mathbf{T}\rightarrow\mathbf{P}\rightarrow M$ be a complete resolution of M. It follows from the
adjointness isomorphism that
\[\begin{array}{rl}
\widehat{\Ext}^i_R(M,N^{\vee})&\cong\hh^i(\Hom_R(\mathbf{T},\Hom_R(N,\E_R(k))))\\
&\cong\hh^i(\Hom_R(\mathbf{T}\otimes_RN,\E_R(k)))\\
&\cong\widehat{\Tor}_i^R(M,N)^{\vee}\\
\end{array}\]
for all $i\in\ZZ$. 
\end{proof}
%%%%%%%%%%%%%%%%%%%%%%%%%%%%%%%%%%%%%%%%%%%%%%%%%%%%%%%%%%%%%%%%%%%%%%%%%%%%%%%%%%%%%%%%%%%%%%%%%%%%%%%%%%%%%%%%%%%%%%%%%%%%%%%%%%%%%%%%%%%%%%
\begin{lem}\label{l3}
Let $M$ and $N$ be $R$--modules such that $\G-dim_R(M)<\infty$. Assume that $x\in R$ is both $R$-regular and $M$-regular and that $xN=0$.
Set $\overline{(-)}:=-\otimes_RR/x$. Then we have the following isomorphisms:
\begin{enumerate}[(i)]
\item{$\widehat{\Ext}^i_R(M,N)\cong\widehat{\Ext}^i_{\overline{R}}(\overline{M},N)$ \text{ for all } $i\in\ZZ$.}
\item{$\widehat{\Tor}_i^R(M,N)\cong\widehat{\Tor}_i^{\overline{R}}(\overline{M},N)$ \text{ for all } $i\in\ZZ$.}
\end{enumerate}
\end{lem}
\begin{proof}
We only prove part (i), the proof of part (ii) is similar. 
Without loss of generality, by Remark \ref{rem}, we may assume that $\G-dim_R(M) = 0$. Let $\mathbf{T}\rightarrow\mathbf{P}\rightarrow M$ be a complete resolution of $M$. It is easy to see that $\overline{\mathbf{T}}\rightarrow\overline{\mathbf{P}}\rightarrow\overline{M}$ is a complete resolution of $\overline{M}$. Therefore, by adjointness isomorphism
\[\begin{array}{rl}
\widehat{\Ext}^i_R(M,N)&\cong\hh^i(\Hom_R(\mathbf{T},N))\\
&\cong\hh^i(\Hom_R(\mathbf{T},\Hom_{\overline{R}}(\overline{R},N))\\
&\cong\hh^i(\Hom_{\overline{R}}(\overline{\mathbf{T}},N))\\
&\cong\widehat{\Ext}^i_{\overline{R}}(\overline{M},N),
\end{array}\]
for all $i\in\ZZ$.
\end{proof}
%%%%%%%%%%%%%%%%%%%%%%%%%%%%%%%%%%%%%%%%%%%%%%%%%%%%%%%%%%%%%%%%
In \cite{Ia}, Iacob proved that the Tate homology is balanced for modules over commutative Noetherian Gorenstein rings. In 
\cite[Theorem 3.7]{CJ1}, Christensen and Jorgensen generalized the Iacob's result for complexes of finite Gorenstein dimension
over any associative ring. In the following we give a short and simple proof of the Christensen and Jorgensen's result for the module case.
\begin{lem}\label{l5}
Let $M$ and $N$ be $R$--modules of finite Gorenstein dimension. Then $\widehat{\Tor}_i^R(M,N)\cong\widehat{\Tor}^R_i(N,M)$ for all $i\in\ZZ$.
\end{lem}
\begin{proof}
Without loss of generality, by Remark \ref{rem}, we may assume that $M$ and $N$ are totally reflexive. Therefore, by (\ref{t}),	
\[\begin{array}{rl}
\widehat{\Tor}_i^R(M,N)&\cong\Tor_i^R(M,N)\\
&\cong\Tor_i^R(N,M)\cong\widehat{\Tor}_i^R(N,M),
\end{array}\]
for all $i>0$. As $\G-dim_R(M)=0=\G-dim_R(N)$, there exist $R$--modules $X_i$ and $Y_i$
of $G$-dimension zero such that $M=\Omega^{-i+1}X_i$ and $N=\Omega^{-i+1}Y_i$ for all $i\leq0$. It
follows from (\ref{t}) that
\[\begin{array}{rl}
\widehat{\Tor}_i^R(M,N)&\cong\Tor_1^R(X_i,N)\cong\Tor_{-i+2}^R(X_i,Y_i)\\
&\cong\Tor_1^R(Y_i,M)\cong\widehat{\Tor}_i^R(N,M),
\end{array}\]
for all $i\leq0$.
\end{proof}
%%%%%%%%%%%%%%%%%%%%%%%%%%%%%%%%%%%%%%%
We are now ready to prove our main Theorem. The proof of the following result was inspired by the proof of Theorem 2.1 in \cite{HJ}.
\begin{thm}\label{t4}
Let $R$ be a Gorenstein local ring of dimension $d$ and let $M$, $N$ be $R$--modules. Assume that $n$ is an integer. The following statements hold:
\begin{enumerate}[(i)]
\item{If $\widehat{\Ext}^i_R(M,N)=0$ for all $i\geq n$, then $\widehat{\Ext}^i_R(N,M)=0$ for all $i<-n$.}
\item{If $\widehat{\Ext}^i_R(M,N)=0$ for all $i<n$, then $\widehat{\Ext}^i_R(N,M)=0$ for all $i\geq d-n$.}
\end{enumerate}
\end{thm}
\begin{proof}
We only prove part (ii), the proof of part (i) is similar. Without loss of generality, by Remark \ref{rem}, we may assume that $M$ and $N$ are totally reflexive. We argue by induction on $d$. If $d=0$, then $\widehat{\Ext}^i_R(M,N)=0$ for all $i<n$ if and only if $\widehat{\Tor}_i^R(M,N^*)=0$ for all $i<n$ by Lemma \ref{l2} and this is equivalent to say that $\widehat{\Ext}^i_R(N,M)=0$ for all $i\geq-n$ by Theorem \ref{tate}(ii) and Lemma \ref{l5}. Assume that $d>0$ and that $x\in R$ a non-zero divisor on $R$, $M$ and $N$. Set $\overline{(-)}:=-\otimes_R\overline{R}$.
Consider the following exact sequence
\begin{equation}\tag{\ref{t4}.3}
0\longrightarrow N\overset{x}{\longrightarrow}N\longrightarrow \overline{N}\longrightarrow0.
\end{equation}
The above exact sequence, induces a doubly infinite long exact sequence
\begin{equation}\tag{\ref{t4}.4}
\cdots\rightarrow\widehat{\Ext}^i_R(M,N)\overset{x}{\rightarrow}\widehat{\Ext}^i_R(M,N)\rightarrow
\widehat{\Ext}^i_R(M,\overline{N})\rightarrow\cdots,
\end{equation}
of Tate cohomology modules by Theorem \ref{tate}(iv). Assume that
$$\widehat{\Ext}^i_R(M,N)=0  \text{ for all } i<n;$$
by (\ref{t4}.4)
$$\Longrightarrow \widehat{\Ext}^i_R(M,\overline{N})=0  \text{ for all } i<n-1;$$
By Lemma \ref{l3},
$$\Longrightarrow \widehat{\Ext}^i_{\overline{R}}(\overline{M},\overline{N})=0  \text{ for all } i<n-1;$$
By the induction hypothesis,
$$\Longrightarrow \widehat{\Ext}^i_{\overline{R}}(\overline{N},\overline{M})=0  \text{ for all } i\geq d-1-n+1;$$
By Lemma \ref{l3},
$$\Longrightarrow \widehat{\Ext}^i_{R}(N,\overline{M})=0  \text{ for all } i\geq d-n;$$
and now from the long exact sequence of $\widehat{\Ext}$ coming from the short exact
sequence $0\rightarrow M\overset{x}{\rightarrow}M\rightarrow\overline{M}\rightarrow0$ and Nakayama's lemma,
$$\Longrightarrow \widehat{\Ext}^i_{R}(N,M)=0  \text{ for all } i\geq d-n.$$
\end{proof}
%%%%%%%%%%%%%%%%%%%%%%%%%%%%%%%
The following is an immediate consequence of Theorem \ref{t4}.
\begin{cor}\label{c3}
Let $R$ be a Gorenstein local ring and let $M$, $N$ be $R$--modules.
The following statements hold:
\begin{enumerate}[(i)]
\item{$\widehat{\Ext}^i_R(M,N)=0$ for all $i\gg0$ if and only if $\widehat{\Ext}^i_R(N,M)=0$ for all $i\ll0$.}
\item{$\widehat{\Ext}^i_R(M,N)=0$ for all $i\in\ZZ$ if and only if $\widehat{\Ext}^i_R(N,M)=0$ for all $i\in\ZZ$.}
\end{enumerate}
\end{cor}
%%%%%%%%%%%%%%%%%%%%%%%%%%%%%%%%%%%%%%%%%%%%%%%%%%%
%%%%%%%%%%%%%%%%%%%%%%%%%%%%%%%%%%%%%%%%%%%%%%%%%%%%%%%%%%%%%%%%%%%%%%%%%%%%%%%%%%%%%%%%%%%%%%%%%%%%%%%%%%%%%%%%%%%%%%%%%%%%%%%%%%%%%
%%%%%%%%%%%%%%%%%%%%%%%%%%%%%%%%%%%%%%%%%%%%%%%%%%%%%%%%%%%%%%%%%%%%%%%%
The following is a generalization of Auslander-Reiten duality Theorem \cite{AA,Bu} (see also \cite[Theorem 3.4]{CT} for a different generalization).
\begin{thm}\label{t1}
Let $(R,\fm,k)$ be a Cohen-Macaulay local ring of dimension $d$ with canonical module $\omega_R$. Assume that $M$, $N$ are $R$--modules of finite Gorenstein dimension and that $\widehat{\Ext}^i_R(M,N)$ has finite length for all $i\in\ZZ$. Then $$\widehat{\Ext}^i_R(M,N)\cong\Hom_R(\widehat{\Ext}^{d-i-1}_R(N,M\otimes_R\omega_R),\E_R(k))$$
for all $i\in\ZZ$, where $\E_R(k)$ is the injective envelope of the residue field of $R$.
\end{thm}
\begin{proof}
Consider the following exact sequences:
\begin{equation}\tag{\ref{t1}.1}
0\rightarrow P\rightarrow X\rightarrow M\rightarrow0,
\end{equation}
\begin{equation}\tag{\ref{t1}.2}
0\rightarrow Q\rightarrow Y\rightarrow N\rightarrow0,
\end{equation}
where $\G-dim_R(X)=0=\G-dim_R(Y)$ and $P$, $Q$ have finite projective dimension.
By Remark \ref{rem}, $\widehat{\Ext}^i_R(M,N)\cong\widehat{\Ext}^i_R(X,Y)$ for all $i\in\ZZ$. As $\G-dim_R(M)<\infty$, $\Tor_i^R(M,\omega_R)=0$ for all $i>0$ by \cite[Proposition 2.5]{Fo}. Therefore, the exact sequence (\ref{t1}.1) induces the following exact sequence
\begin{equation}\tag{\ref{t1}.3}
0\rightarrow P\otimes_R\omega_R\rightarrow X\otimes_R\omega_R\rightarrow M\otimes_R\omega_R\rightarrow0,
\end{equation}
Note that $\id_R(P\otimes_R\omega_R)<\infty$ and so  $\widehat{\Ext}^i_R(N,P\otimes_R\omega_R)=0$ for all $i\in\ZZ$ by Theorem \ref{tate}(iii).
Therefore the excat sequences (\ref{t1}.2) and (\ref{t1}.3) induce the following isomorphisms $\widehat{\Ext}^i_R(N,M\otimes_R\omega_R)\cong\widehat{\Ext}^i_R(N,X\otimes_R\omega_R)\cong\widehat{\Ext}^i_R(Y,X\otimes_R\omega_R)$ for all $i\in\ZZ$ by Theorem \ref{tate}(iv) ,(v).
So we may assume that $M$ and $N$ are maximal Cohen-Macaulay (equivalentely, totally reflexive). Note that $\G-dim_R(\Omega^{i}M)=0$ for all $i\in\ZZ$ and so
$\G-dim_R(\Tr\Omega^{i}M)=0$ for all $i\in\ZZ$ by Theorem \ref{G}(i).
It follows from (\ref{t}) and Theorem \ref{tate}(ii) that
\[\begin{array}{rl}\tag{\ref{t1}.4}
\Tor_j^R(\Tr\Omega^{i}M,N)&\cong\widehat{\Tor}_j^R(\Tr\Omega^{i}M,N)\\
&\cong\widehat{\Tor}_{j-2}^R((\Omega^{i}M)^*,N)\\
&\cong\widehat{\Ext}^{-j+1}_R(\Omega^{i}M,N)\cong
\widehat{\Ext}^{i-j+1}_R(M,N),
\end{array}\]
for all $j>0$ and $i\in\ZZ$. 
Hence, by (\ref{t1}.4) and our assumption, $\Tor_j^R(\Tr\Omega^{i}M,N)$ has finite length for all $i\in\ZZ$ and $j>0$.
By \cite[Theorem 10.62]{R}, there is a third quadrant spectral sequence:
$$\E^{p,q}_2=\Ext^p_R(\Tor_q^R(\Tr\Omega^{i}M,N),\omega_R)\Rightarrow\Ext^{p+q}_R(N,(\Tr\Omega^{i}M)^\dag),$$
where $(-)^\dag:=\Hom_R(-,\omega_R)$. As $\Tor_q^R(\Tr\Omega^{i}M,N)$ has finite length for all $q>0$ and $i\in\ZZ$,
$\E^{p,q}_2=0$ if $p\neq d$. Hence the spectral sequence collapses and so
\begin{equation}\tag{\ref{t1}.6}
\Ext^d_R(\Tor_j^R(\Tr\Omega^{i}M,N),\omega_R)\cong\Ext^{d+j}_R(N,(\Tr\Omega^{i}M)^\dagger)
\end{equation}
for all $j>0$ and $i\in\ZZ$. As $\Tr\Omega^iM$ is totally reflexive, the exact sequence (\ref{1.1}),
$$0\rightarrow(\Omega^iM)^*\rightarrow P_0^*\rightarrow P_1^*\rightarrow\Tr\Omega^iM\rightarrow0,$$
induces the following exact sequence:
$$0\rightarrow(\Tr\Omega^iM)^{\dagger}\rightarrow P_1^{*\dagger}\rightarrow P_0^{*\dagger}\rightarrow(\Omega^iM)^{*\dagger}\rightarrow0.$$
From the above exact sequence we get the following isomorphism:
\begin{equation}\tag{\ref{t1}.7}
\widehat{\Ext}^{j+2}_R(N,(\Tr\Omega^{i}M)^\dagger)\cong\widehat{\Ext}^{j}_R(N,(\Omega^iM)^{*\dagger}),
\end{equation}
for all $i\in\ZZ$ and $j\in\ZZ$, by Theorem \ref{tate}(iii), (v). As $\G-dim_R(\Omega^iM)=0$, $\Omega^iM\otimes_R\omega_R$ is maximal Cohen-Macaulay 
by \cite[Theorem 2.13]{Y} and so
\begin{equation}\tag{\ref{t1}.8}
\Omega^iM\otimes_R\omega_R\cong(\Omega^iM\otimes_R\omega_R)^{\dagger\dagger}\cong (\Omega^iM)^{*\dagger}\text{ for all } i\in\ZZ.
\end{equation}
Note that $\Tor_j^R(\Omega^iM,\omega_R)=0$ for all $i\in\ZZ$ and $j>0$ by \cite[Proposition 2.5]{Fo}. Hence, the exact sequence $0\rightarrow\Omega^iM\rightarrow F\rightarrow\Omega^{i-1}M\rightarrow0$, induces the  following exact sequence
$$0\rightarrow\Omega^iM\otimes_R\omega_R\rightarrow F\otimes_R\omega_R\rightarrow\Omega^{i-1}M\otimes_R\omega_R\rightarrow0.$$
From the above exact sequence we get the following isomorphism:
\begin{equation}\tag{\ref{t1}.9}
\widehat{\Ext}^{j-1}_R(N,\Omega^{i-1}M\otimes_R\omega_R)\cong\widehat{\Ext}^{j}_R(N,\Omega^{i}M\otimes_R\omega_R),
\end{equation}
for all $i\in\ZZ$ and $j\in\ZZ$, by Theorem \ref{tate}(iii), (iv).
As $\Tor_1^R(\Tr\Omega^{i}M,N)$ has finite length, we get the following isomorphisms by (\ref{t1}.4)
\begin{equation}\tag{\ref{t1}.10}
\widehat{\Ext}^{i}_R(M,N)\cong\Tor_1^R(\Tr\Omega^{i}M,N)
\cong\Gamma_{\fm}(\Tor_1^R(\Tr\Omega^{i}M,N)).
\end{equation}
Now the assertion follows from the Local duality Theorem,  (\ref{t1}.6), (\ref{t1}.7), (\ref{t1}.8), (\ref{t1}.9) and (\ref{t1}.10).

\end{proof}
%%%%%%%%%%%%%%%%%%%%%%%%%%%%%%%%%%%%%%%%%%%%%%%%%%%%%%%%%%%%%%%%%%%%%%%%%%%%%%%%%%%%%%%%%%
%%%%%%%%%%%%%%%%%%%%%%%%%%%%%%%%%%%%%%%%%%%%%%%%%%%%%%%%%%%%%%%%%%%%%%%%%%%%%%%%%%%%%%%%%%%%%%%
Let $(R,\fm,k)$ be a Gorenstein local ring of dimension $d$ and let $M$ be an $R$--module of finite projective dimension. In \cite[Corollary 3.2]{Fox}, Foxby proved that $\beta_i^R(M)=\mu^{d-i}_R(M)$ for all $i$. The $i$th stable Betti number and $i$th stable Bass number was introduced by Avramov and Martsinkovsky \cite{AM}. For an $R$--module $N$,
$$\widehat{\beta}_i^R(N)=\rank_k(\widehat{\Ext}^i_R(N,k)) \text{ and } \widehat{\mu}^i_R(N)=\rank_k(\widehat{\Ext}^i_R(k,N)).$$
Note that $\widehat{\beta}_i^R(N)=\beta_i^R(N)$ for all $i>\G-dim_R(N)$ and also $\widehat{\mu}^i_R(N)=\mu^i_R(N)$ for all $i>d$.
Let $N$ be an $R$--module of infinite projective dimension. In \cite[Theorem 10.3]{AM}, it is shown that $\widehat{\mu}^i_R(N)=\beta_{d-i-1}^R(N)$ for all $i<\depth_R(N)-1$. The following is an immediate consequence of Theorem \ref{t1}.
\begin{cor}\label{c1}
Let $R$ be a Gorenstein local ring of dimension $d$ and let $M$ be an $R$--module. Then $\widehat{\beta}_i^R(M)=\widehat{\mu}^{d-i-1}_R(M)$ for all $i\in\ZZ$.
\end{cor}
The following is an immediate consequence of Corollary \ref{c1} and \cite[Theorem 9.1]{AM}.
\begin{cor}\label{c6}
Let $R$ be a Gorenstein local ring and let $\fa$ be a Gorenstein
ideal of $R$. Then $\widehat{\beta}_i^R(R/\fa)=\widehat{\mu}^{n+i}_R(R/\fa)$ for all $i\in\ZZ$, where $n=\dim_R(R/\fa)$.
In particular, $\beta_i^R(R/\fa)=\mu^{n+i}_R(R/\fa)$ for all $i>\G-dim_R(R/\fa)$.
\end{cor}
%%%%%%%%%
%%%%%%%%%%%%%%%%%%%%%%%%%%%%%%%%%%%%%%%%%%%%%%%%%%%%%%%%%%%%%%%%%%
%%%%%%%%%%%%%%%%%%%%%%%%%%%%%%%%%%%%%%%%%%%%%%%%%%%%%%%%%%%%%%%%%%%%%%%%%%%%%%%
\section{Applications of main result}
%%%%%%%%%%%%%%%%%%%%%%%%%%%%%%%%%%%%%%%%%%%%%%%%%%%%%%%%%%%%%%%%%%%%%%%%%%%%%%%%%
In \cite[Thoerem 4.1]{J}, J${\o}$rgensen proved that the symmetry in the vanishing of Ext holds for modules of finite complete intersection dimension over local Gorenstein rings. Bergh generalized the J${\o}$rgensen's result for modules with reducible complexity (see \cite[Theorem 3.5]{B2}).
In this section, it is shown that the symmetry in the vanishing of Ext holds for modules over local Gorenstein rings, provided that one of them has reducible complexity.
\begin{prop}\label{l4}
Let $M$ and $N$ be $R$--modules such that $\G-dim_R(M)<\infty$. Assume that either $M$ or $N$ has reducible complexity. Then the following are equivalent:
\begin{enumerate}[(i)]
\item{$\widehat{\Ext}^i_R(M,N)=0$ for all $i\gg0$;}
\item{$\widehat{\Ext}^i_R(M,N)=0$ for all $i\ll0$;}
\item{$\widehat{\Ext}^i_R(M,N)=0$ for all $i\in\ZZ$.}
\end{enumerate}
\end{prop}
\begin{proof}
We only prove the implication (ii)$\Rightarrow$(iii), when $N$ has reducible complexity. The proof of other implications are similar.
Set $c=\cx_R(N)$. We argue by induction on $c$. If $c=0$ then $\pd_R(N)<\infty$ and the assertion follows from Theorem \ref{tate}(i). Now let $c>0$ and $\eta\in\Ext^*_R(N,N)$ reduces the complexity of $N$. Consider the exact sequence
\begin{equation}\tag{\ref{l4}.1}
0\rightarrow N\rightarrow K_{\eta}\rightarrow\Omega^q_R(N)\rightarrow0,
\end{equation}
where $q=|\eta|-1$ and $\cx_R(K_{\eta})<c$. The exact sequence (\ref{l4}.1), induces a doubly infinite long exact sequence
\begin{equation}\tag{\ref{l4}.2}
\cdots\rightarrow\widehat{\Ext}^i_R(M,N)\rightarrow\widehat{\Ext}^i_R(M,K_{\eta})\rightarrow\widehat{\Ext}^i_R(M,\Omega^q_R(N))\rightarrow\cdots
\end{equation}
of Tate cohomology modules by Theorem \ref{tate}(iv).
Note that
\begin{equation}\tag{\ref{l4}.3}
\widehat{\Ext}^i_R(M,\Omega^q_R(N))\cong\widehat{\Ext}^{i-q}_R(M,N) \text{ for all } i\in\ZZ.
\end{equation}
It follows from (\ref{l4}.2), (\ref{l4}.3) and (ii) that $\widehat{\Ext}^i_R(M,K_{\eta})=0$ for all $i\ll0$. By induction hypothesis,
$\widehat{\Ext}^i_R(M,K_{\eta})=0$ for all $i\in\ZZ$. Therefore,
\begin{equation}\tag{\ref{l4}.4}
\widehat{\Ext}^i_R(M,N)\cong\widehat{\Ext}^{i-q-1}_R(M,N) \text{ for all } i\in\ZZ.
\end{equation}
As $\widehat{\Ext}^i_R(M,N)=0$ for all $i\ll0$, the assertion is clear by (\ref{l4}.4).
\end{proof}
%%%%%%%%%%%%%%%%%%%%%%%%%%%%%%%%%%%%%%%%%%%%%%%%%%%%%%%%%%%%%%%%%%%%%
\begin{cor}\label{c7}
Let $R$ be a Gorenstein local ring and let $M$, $N$ be $R$--modules. Assume that $N$ has reducible complexity$($e.g. $\CI_R(N)<\infty)$. Then the following are equivalent:
\begin{enumerate}[(i)]
\item{$\Ext^i_R(M,N)=0$ for $i\gg0$.}
\item{$\Ext^i_R(N,M)=0$ for $i\gg0$.}
\end{enumerate}
\end{cor}
\begin{proof}
By Proposition \ref{l4}, $\Ext^i_R(M,N)=0$ for all $i\gg0$,
$$\Longleftrightarrow \widehat{\Ext}^i_R(M,N)=0 \text{ for all } i\in\ZZ,$$
$$\Longleftrightarrow \widehat{\Ext}^i_R(N,M)=0 \text{ for all } i\in\ZZ, \text{ by Corollary \ref{c3}}$$
$$\Longleftrightarrow \Ext^i_R(N,M)=0 \text{ for all } i\gg0, \text{ by Proposition \ref{l4}}.$$
\end{proof}

%%%%%%%%%%%%%%%%%%%%%%%%%%%%%%%%%%%%%%%%%%%%%%%%%%%%%%%%%%%%%%%%%%%%%%%%%%%%%%
%%%%%%%%%%%%%%%%%%%%%%%%%%%%%%%%%%%%%%%%%%%%%%%%%%%%%%%%%%%%%%%%%%%%%%%%%%%%%%%%%%%%%%%%%
Let $R$ be a ring and let $M$ and $N$ be $R$--modules. Assume that $M$ has finite Gorenstein dimension. In \cite[Theorem 6.1]{CJ},
Christensen and Jorgensen proved that
$$\sup\{i\mid\Ext^i_R(M,N)\neq0\}=\depth R-\depth_R(M),$$
provided that $\widehat{\Ext}^i_R(M,N)=0$ for all $i\in\ZZ$.
In the following, we prove the same conclusion under
weaker hypotheses:
%%%%%%%%%%%%%%%%%%%%%%%%%%%%%%%%%%%%%%%%%%%%%%%%%%%%%%%%%%%%%%%%%%%%%%%%%%%%%%%%%%%%%%%
\begin{lem}\label{l1}
Let $M$ and $N$ be nonzero $R$--modules such that
$\G-dim_R(M)<\infty$. If $\widehat{\Ext}^i_R(M,N)=0$ for all
$i\geq\G-dim_R(M)-1$, then the following statements hold.
\begin{enumerate}[(i)]
\item{$\sup\{i\mid\Ext^i_R(M,N)\neq0\}=\depth R-\depth_R(M)$;}
\item{$\Ext^{\tiny{\G-dim_R(M)}}_R(M,N)\cong\Ext^{\tiny{\G-dim_R(M)}}_R(M,R)\otimes_RN$.}
\end{enumerate}
\end{lem}
\begin{proof}
Set $n=\G-dim_R(M)$. Note that
$\Ext^i_R(M,N)\cong\widehat{\Ext}^i_R(M,N)$ for all $i>n$. Therefore
$\Ext^i_R(M,N)=0$ for all $i>n$. Set $L=\Omega^n M$. Note that $L^*\approx\Omega^2\Tr L$. Hence
\[\begin{array}{rl}\tag{\ref{l1}.1}
\widehat{\Tor}_i^R(\Tr L,N)&\cong\widehat{\Tor}_{i-2}^R(L^*,N)\\
&\cong\widehat{\Ext}^{-i+1}_R(L,N)\\
&\cong\widehat{\Ext}^{-i+1+n}_R(M,N),
\end{array}\]
for all $i\in\ZZ$, by Theorem \ref{tate}(ii). As $\G-dim_R(L)=0$, by Theorem \ref{G}(i), $\G-dim_R(\Tr L)=0$ and so by (\ref{t})
\begin{equation}\tag{\ref{l1}.2}
\Tor_i^R(\Tr L,N)\cong\widehat{\Tor}_i^R(\Tr L,N)  \text{ for } i>0.
\end{equation}
It follows from our assumption, (\ref{l1}.1) and (\ref{l1}.2) that $\Tor_i^R(\Tr L,N)=0$ for $i=1,2$.
Therefore, $\Ext^n_R(M,R)\otimes_RN\cong\Ext^n_R(M,N)$ by Theorem \ref{a1} and so the assertion is clear by Theorem \ref{G}(ii).
\end{proof}
Note that the vanishing of negative and positive Tate cohomology
are two distinct conditions in general. By \cite[Theorem 4.1]{JS1}, there
exists an artinian Gorenstein local ring $R$ and finitely generated $R$--modules $M$ and
$N$ such that $\widehat{\Ext}^i_R(M,N)=0$ for all $i>0$ and $\widehat{\Ext}^i_R(M,N)\neq0$ for all $i<0$.
%%%%%%%%%%%%%%%%%%%%%%%%%%%%%%%%%%%%%%%%%%%%%%%%%%%%%%%%%%%%%%%%%%%%%%%%%%%%%%%%%%%%%%%%%%%
The following is an immediate consequence of Theorem \ref{t4} and Lemma \ref{l1}.
\begin{cor}\label{c4}
Let $R$ be a Gorenstein local ring and let $M$, $N$ be non-zero $R$--modules.
Assume that $\widehat{\Ext}^i_R(M,N)=0$ for all $i\leq\depth_R(N)$. Then $$\sup\{i\mid\Ext^i_R(N,M)\neq0\}=\depth R-\depth_R(N).$$
\end{cor}
%%%%%%%%%%%%%%%%%%%%%%%%%%%%%%%%%%%%%%%%%%%%%%%%%%%%%%%%%%%%%%%%%%%%%%%%%%%%%%%%%%%%%%%%
\section{G-perfect modules}
Recall that an $R$--module $M$ is $G$-perfect if
$\gr_R(M)=\G-dim_R(M)$. Over a Gorenstein ring,
these modules are precisely the Cohen-Macaulay modules.
For a $G$-perfect $R$--module $M$ of
Gorenstein dimension $n$, we set $M^{\dagger}:=\Ext^{n}_R(M,R)$.
In the following, we collect some basic properties of $G$-perfect modules (see \cite{G} for more details).
\begin{thm}\label{perfect}
Let $M$ be a $G$-perfect $R$--module of Gorenstein dimension $n$.
The following statements hold.
\begin{enumerate}[(i)]
\item{$M^{\dagger}$ is a $G$-perfect $R$--module of Gorenstein dimension $n$;}
\item{$M\cong M^{\dagger\dagger}$;}
\item{$\ann_R(M)=\ann_R(M^{\dagger}).$}
\end{enumerate}
\end{thm}
Let $M$ be a perfect $R$--module and $N$ an $R$--module. The following equality is well-known (see for example \cite[Proposition 2.9.3]{H}),
$$\sup\{i\geq0\mid\Tor_i^R(M,N)\neq0\}+\inf\{i\geq0\mid\Ext^i_R(M,N)\neq0\}=\pd_R(M).$$
In the following we generalize the above equality for $G$-perfect modules.
%%%%%%%%%%%%%%%%%%%%%%%%%%%%%%%%%%%%%%%%%%%%%%%%%%%%%%%%%%%%%%%%%%%%%%%%%%%%%%%%%%%%%%
Recall that $\gr_R(M,N)=\inf\{i\geq0\mid\Ext^i_R(M,N)\neq0\}$.
\begin{thm}\label{t3}
Let $M$ and $N$ be non-zero $R$--modules. Assume that $M$ is $G$-perfect.
If $\widehat{\Tor}_i^R(M,N)=0$ for all $i\geq0$(e.g., either $\pd_R(M)<\infty$ or $\id_R(N)<\infty$), then
$$\sup\{i\geq0\mid\Tor_i^R(M,N)\neq0\}+\inf\{i\geq0\mid\Ext^i_R(M,N)\neq0\}=\G-dim_R(M).$$
\end{thm}
\begin{proof}
Set $n=\G-dim_R(M)$ and $t=\sup\{i\mid\Tor_i^R(M,N)\neq0\}$. Since
$\Tor_i^R(M,N)\cong\widehat{\Tor}_i^R(M,N)$ for all $i>n$, it
follows that $t\leq n$. Set $M^{\dag}:=\Ext^n_R(M,R)$.
By \cite[Proposition 6]{F}, there exists the following exact sequence:
\[\begin{array}{lr}\tag{\ref{t3}.1}
0\longrightarrow\Ext^1_R(\Tr\Omega^{n}M,N)\longrightarrow\Tor_{n}^R(M,N)\longrightarrow\Hom_R(M^{\dag},N)\longrightarrow\cdots\\
\cdots\longrightarrow\Ext^{j}_R(\Tr\Omega^{n}M,N)\longrightarrow\Tor_{n+1-j}^R(M,N)\longrightarrow\Ext^{j-1}_R(M^{\dag},N)\longrightarrow\cdots\\
\cdots\rightarrow\Ext^{n+1}_R(\Tr\Omega^{n}M,N)\rightarrow M\otimes_RN\rightarrow\Ext^{n}_R(M^{\dag},N)\rightarrow\Ext^{n+2}_R(\Tr\Omega^{n}M,N)\rightarrow0.
\end{array}\]
Note that $\G-dim_R(\Omega^nM)=0$ and $\Omega^2\Tr\Omega^{n}M\approx(\Omega^{n}M)^*$. Hence, by 
Theorem \ref{tate}(ii)
\[\begin{array}{rl}\tag{\ref{t3}.2}
\widehat{\Tor}_i^R(M,N)&\cong\widehat{\Tor}_{i-n}^R(\Omega^nM,N)\\
&\cong\widehat{\Ext}^{-i+n-1}_R((\Omega^{n}M)^*,N)\\
&\cong\widehat{\Ext}^{-i+n+1}_R(\Tr\Omega^{n}M,N).
\end{array}\]
for all $i\in\ZZ$. It follows from (\ref{t3}.2) and our assumption that
\begin{equation}\tag{\ref{t3}.3}
\Ext^i_R(\Tr\Omega^{n}M,N)\cong\widehat{\Ext}^i_R(\Tr\Omega^{n}M,N)=0  \text{ for all } 1\leq i\leq n+1.
\end{equation}
Therefore, by (\ref{t3}.1) and (\ref{t3}.3),
\begin{equation}\tag{\ref{t3}.4}
\Tor_{n-i}^R(M,N)\cong\Ext^{i}_R(M^{\dag},N) \text{ for all } 0\leq i\leq n-1,
\end{equation}
and also we have the following exact sequence
\begin{equation}\tag{\ref{t3}.5}
0\longrightarrow M\otimes_RN\longrightarrow\Ext^{n}_R(M^{\dag},N)
\end{equation}
It follows from (\ref{t3}.5) that $\gr_R(M^{\dagger},N)\leq n$.
Therefore, $t=n-\gr_R(M^{\dag},N)$ by (\ref{t3}.4). On the other hand, $M$ is
$G$-perfect and so $\ann_R(M^{\dagger})=\ann_R(M)$ by Theorem \ref{perfect}. Therefore
\[\begin{array}{rl}
\gr_R(M^{\dagger},N)&=\gr_R(\ann_R(M^{\dagger}),N))\\
&=\gr_R(\ann_R(M),N)\\
&=\gr_R(M,N).
\end{array}\]
\end{proof}
%%%%%%%%%%%%%%%%%%%%%%%%%%%%%%%%%%%%%%%%%%%%%%%%%%%%%%%%%%%%%%%%%%%%%%%%%%%%%%%%%%%%%%%%%%
The following is an immediate consequence of Theorem \ref{t3}, \cite[Proposition 3.2]{CJ} and \cite[Theorem 4.9]{AvBu}.
\begin{cor}
Let $M$ be a $G$-perfect $R$--module and $N$ an $R$--module.
Assume that $\Tor_i^R(M,N)=0$ for $i\gg0$ and that one of the
following conditions hold:
\begin{enumerate}[(i)]
\item{$R$ is AB ring;}
\item{Either $M$ or $N$ has finite complete intersection dimension;}
\end{enumerate}
Then,
$\sup\{i\geq0\mid\Tor_i^R(M,N)\neq0\}+\inf\{i\geq0\mid\Ext^i_R(M,N)\neq0\}=\G-dim_R(M).$
\end{cor}
%%%%%%%%%%%%%%%%%%%%%%%%%%%%%%%%%%%%%%%%%%%%%%%%%%%%%%%%%%%%%%%%%%%%%%%%%%%%%%%%%%%%%%%%%
%%%%%%%%%%%%%%%%%%%%%%%%%%%%%%%%%%%%%%%%%%%%%%%%%%%%%%%%%%%%%%%%%%%%%%%
%%%%%%%%%%%%%%%%%%%%%%%%%%%%%%%%%%%%%%%%%%%%%%%%%%%%%%%%%%%%%%%%%%%%%%%%%%%%%%%%%%%%%%%%%%%%
%%%%%%%%%%%%%%%%%%%%%%%%%%%%%%%%%%%%%%%%%%%%%%%%%%%%%%%%%%%%%%%%%%%%%
The following result plays a crucial role in this paper.
\begin{prop}\label{pr1}
Let $M$ be a $G$-perfect $R$--module of Gorenstein dimension $n$ and
let $N$ be an $R$--module. Then we have the following isomorphisms.
\begin{enumerate}[(i)]
\item{$\widehat{\Ext}^i_R(M^{\dagger},N)\cong\widehat{\Tor}_{-i-1+n}^R(M,N)$,}
\item{$\widehat{\Ext}^i_R(M,N)\cong\widehat{\Tor}_{-i-1+n}^R(M^{\dagger},N)$,}
\end{enumerate}
for all $i\in\mathbb{Z}$.
\end{prop}
\begin{proof}
If $n=0$ then the assertion follows from Theorem \ref{tate}(ii). So assume $n>0$ and consider the following exact sequence,
\begin{equation}\tag{\ref{pr1}.1}
0\rightarrow M^{\dagger}\rightarrow\Tr\Omega^{n-1}M\rightarrow
L\rightarrow0,
\end{equation}
where $L\approx\Omega\Tr\Omega^{n}M$. As $\gr_R(M)=n$, it is easy to see that
$\pd_R(\Tr\Omega^{n-1}M)<\infty$. Therefore,
$\widehat{\Ext}^i_R(\Tr\Omega^{n-1}M,N)=0$ for all $i\in\mathbb{Z}$ by Theorem \ref{tate}(i).
From the exact sequence (\ref{pr1}.1) and Theorem \ref{tate}(v),
we obtain the following isomorphism
\begin{equation}\tag{\ref{pr1}.2}
\widehat{\Ext}^i_R(M^{\dagger},N)\cong\widehat{\Ext}^{i+1}_R(\Omega\Tr\Omega^{n}M,N)\cong\widehat{\Ext}^{i+2}_R(\Tr\Omega^{n}M,N),
\end{equation}
for all $i\in\mathbb{Z}$. Note that
$\G-dim_R(\Tr\Omega^{n}M)=0$ by Theorem \ref{G}(i). It follows from Theorem \ref{tate}(ii) that
\begin{equation}\tag{\ref{pr1}.3}
\widehat{\Ext}^i_R(\Tr\Omega^{n}M,N)\cong\widehat{\Tor}_{-i-1}^R((\Tr\Omega^{n}M)^*,N),
\end{equation}
for all $i\in\mathbb{Z}$. As
$(\Tr\Omega^{n}M)^*\approx\Omega^2\Tr\Tr\Omega^{n}M\approx\Omega^{n+2}M$,
it follows from (\ref{pr1}.3) that
\begin{equation}\tag{\ref{pr1}.4}
\widehat{\Ext}^i_R(\Tr\Omega^{n}M,N)\cong\widehat{\Tor}_{-i+1+n}^R(M,N)
\end{equation}
Now (i) is clear by (\ref{pr1}.2) and (\ref{pr1}.4).

By Theorem \ref{perfect}, $M^{\dagger}$ is $G$-perfect of Gorenstein
dimension $n$ and $M\cong M^{\dagger\dagger}$. Now (ii) follows from
(i) by replacing $M$ by $M^{\dagger}$.
\end{proof}
%%%%%%%%%%%%%%%%%%%%%%%%%%%%%%%%%%%%%%%%%%%%%%%%%%%%%%%%%%%%%%%%%%%%%%%%%
Let $R$ be a Gorenstein local ring. Recall, an ideal $I$ in $R$ is said to be a
Gorenstein ideal if $R/I$ is a Gorenstein local ring. The following is an immediate consequence of Proposition \ref{pr1} and Corollary \ref{c3}.
\begin{prop}\label{pr2}
Let $R$ be a Gorenstein ring and let $\fa$ be a Gorenstein ideal of $R$. Assume that $M$ is an $R$--module.
Then $\Ext^i_R(M,R/\fa)=0$ for all $i\gg0$ if and only if $\Tor_i^R(M,R/\fa)=0$ for all $i\gg0$.
\end{prop}
%%%%%%%%%%%%%%%%%%%%%%%%%%%%%%%%%%%%%%%%%%%%%%%%%%%%%%%%%%%%%%%%%%%%%%%%%%%%%%%%%%
\begin{cor}\label{c2}
Let $R$ be a Gorenstein ring and let $\fa$, $\fb$ be Gorenstein ideals of $R$. The following conditions are equivalent:
\begin{enumerate}[(i)]
\item{$\Ext^i_R(R/\fa,R/\fb)=0$ for $i\gg0$;}
\item{$\Ext^i_R(R/\fb,R/\fa)=0$ for $i\gg0$;}
\item{$\Tor_i^R(R/\fa,R/\fb)=0$ for $i\gg0$;}
\end{enumerate}
\end{cor}
\begin{proof}
(i)$\Leftrightarrow$(ii). Set $n=\G-dim(R/\fa)$, $m=\G-dim(R/\fb)$. By Proposition \ref{pr1}, we have
\begin{equation}\tag{\ref{c2}.1}
\widehat{\Ext}^i_R(R/\fa,R/\fb)\cong\widehat{\Tor}_{-i-1+n}^R(R/\fa,R/\fb)\cong\widehat{\Ext}^{i+m-n}_R(R/\fb,R/\fa),
\end{equation}
for all $i\in\ZZ$. Now the equivalence of (i), (ii) is clear by (\ref{c2}.1).

(ii)$\Leftrightarrow$(iii). Follows from Proposition \ref{pr2}.
\end{proof}
%%%%%%%%%%%%%%%%%%%%%%%%%%%%%%%%%%%%%%%%%%%%%%%%%%%%%%%%%%%%%%%%%%%%%%%%%%%%%%%%%
\begin{cor}
Let $R$ be a Gorenstein local ring and let $\fa$ be a Gorenstein
ideal of $R$. Assume that either $\CI_R(R/\fa)<\infty$ or $R$ is AB
ring. Then $$\pd_R(R/\fa)=\sup\{i\mid\Tor_i^R(R/\fa,R/\fa)\neq0\}.$$
\end{cor}
\begin{proof}
Assume that $\Tor_i^R(R/\fa,R/\fa)=0$ for all $i\gg0$. It follows from Corollary \ref{c2} that
$\Ext^i_R(R/\fa,R/\fa)=0$ for all $i\gg0$ and so $\pd_R(R/\fa)<\infty$ by
\cite[Corollary 4.4]{AY} and \cite[Proposition 3.2]{CJ}. Set $n=\pd_R(R/\fa)$. By Theorem \ref{a1},
there exists the following exact sequence:
\[\begin{array}{lr}
0\longrightarrow\Ext^1_R(\Tr\Omega^{n}(R/\fa),R/\fa)\longrightarrow\Tor_{n}^R(R/\fa,R/\fa)\longrightarrow\\
\longrightarrow\Hom_R(\Ext^n_R(R/\fa,R),R/\fa)\longrightarrow\Ext^{2}_R(\Tr\Omega^{n}(R/\fa),R/\fa).
\end{array}\]
As $\pd_R(R/\fa)=n$, $\Tr\Omega^{n}(R/\fa)\approx0$ and so $\Tor_{n}^R(R/\fa,R/\fa)\cong\Hom_R(\Ext^n_R(R/\fa,R),R/\fa)$.
Since $R/\fa$ is Gorenstein, $\Ext^n_R(R/\fa,R)\cong R/\fa$. Therefore,  $\Tor_{n}^R(R/\fa,R/\fa)\neq0$.
\end{proof}
%%%%%%%%%%%%%%%%%%%%%%%%%%%%%%%%%%%%%%%%%%%%%%%%%%%%%%%%%%%%%%%%%%%%%%%%%%%%%%

%%%%%%%%%%%%%%%%%%%%%%%%%%%%%%%%%%%%%%%%%%%%%%%%%%%%%%%%%%%%%%%%%%%%%%%%%%%%%%%%
\section{Vanishing of Tate (co)homology for linked modules}
%%%%%%%%%%%%%%%%%%%%%%%%%%%%%%%%%%%%%%%%%%%%%%%%%%%%%%%%%%%%%%%%%%%%%%%%%%%%%%%%%%
The theory of linkage of algebraic varieties introduced by Peskine
and Szpiro \cite{PS}. Recall that two ideals $\fa$ and $\fb$ in a
Gorenstein local ring $R$ are said to be linked by a Gorenstein ideal $\fc$ if $\fc\subseteq\fa\cap\fb$,
$\fa=(\fc:\fb)$ and $\fb=(\fc:\fa)$.
Martsinkovsky and Strooker \cite{MS} generalized
the notion of linkage for modules over non-commutative semiperfect
Noetherian rings by introducing the operator $\lambda=\Omega\Tr$.
They show that ideals $\fa$ and $\fb$ are linked by zero ideal if
and only if $R/\fa$ and $R/\fb$ are related to each other through
the operator $\lambda$; more precisely, $R/\fa\cong\lambda (R/\fb)$
and $R/\fb\cong\lambda (R/\fa)$ \cite[Proposition 1]{MS}.
In \cite{Pu}, Puthenpurakal studied the vanishing of (co)homology for linked modules. In this section, we study the vanishing of Tate (co)homology for linked modules over Gorenstein rings.
As a consequence of our main result, we can generalize some results of Puthenpurakal.

First, we recall the definition of linkage of module.
\begin{dfn}\label{def1}
\cite[Definition 4]{MS} An $R$--module
$M$ is said to be \emph{linked} to an $R$--module $N$, by an ideal $\fc$ of $R$, if $\fc\subseteq\ann_R(M)\cap\ann_R(N)$, $M\cong\lambda_{R/\fc}N$
and $N\cong\lambda_{R/\fc}M$. In this situation we denote $M\underset{\fc}{\thicksim}N$. An $R$--module
$M$ is called \emph{self-linked} if $M$ is linked to itself.
\end{dfn}
Let $M$ be an $R$--module which is linked by an ideal $\fc$ of $R$. By \cite[Proposition 3]{MS}, $M$ is stable as an $R/\fc$--module. Let $P_1\rightarrow P_0\rightarrow M\rightarrow0$ be a minimal projective presentation of $M$ over $R/\fc$. Then $\Hom_{R/\fc}(P_0,R/\fc)\rightarrow \Hom_{R/\fc}(P_1,R/\fc)\rightarrow\Tr_{R/\fc}M\rightarrow0$ is a minimal projective presentation of $\Tr_{R/\fc}M$ (see \cite[Theorem 32.13]{AF}). Therefore we obtain the following exact sequence:
\begin{equation}\label{ttt}
0\rightarrow\Hom_{R/\fc}(M,R/\fc)\rightarrow\Hom_{R/\fc}(P_0,R/\fc)\rightarrow\lambda_{R/\fc}M\rightarrow0.
\end{equation}

Let $\cm^n(R)$ be the full subcategory of Cohen-Macaulay $R$--modules of codimension $n$. Note that for an $R$--module $M$ over Gorenstein local ring $R$ we have $M\in\cm^n(R)$ if and only if $M$ is $G$-perfect of Gorenstein dimension $n$. For an $R$--module $M\in\cm^n(R)$, set $M^\dagger:=\Ext^n_R(M,R)$.

The following is a generalization of \cite[Theorem 1.2]{Pu}.
\begin{prop}\label{t44}
Let $R$ be a Gorenstein ring, $\fc$ a perfect Gorenstein ideal of $R$ and $L$ an $R$--module. Assume that $M\in\cm^n(R)$ and that $M\underset{\fc}{\sim}N$.
The following statements hold:
\begin{enumerate}[(i)]
\item{$\widehat{\Ext}^i_R(L,M)=0$ for $i\gg0$ if and only if $\widehat{\Tor}_i^R(L,N)=0$ for $i\gg0$.}
\item{$\widehat{\Ext}^i_R(L,M)=0$ for $i\ll0$ if and only if $\widehat{\Tor}_i^R(L,N)=0$ for $i\ll0$.}
\item{$\widehat{\Ext}^i_R(L,M)=0$ for $i\in\ZZ$ if and only if $\widehat{\Tor}_i^R(L,N)=0$ for $i\in\ZZ$.}
\end{enumerate}
\end{prop}
\begin{proof}
We only prove part (ii). The other parts can be proved similarly. By Corollary \ref{c3}, Proposition \ref{pr1} and Lemma \ref{l5},
\[\begin{array}{rl}\tag{\ref{t44}.1}
\widehat{\Ext}^i_R(L,M)&=0 \text{ for all } i\ll0\\
\Longleftrightarrow\widehat{\Ext}^{i}_R(M,L)&=0 \text{ for all } i\gg0\\
\Longleftrightarrow\widehat{\Tor}_{i}^R(M^{\dagger},L)&=0  \text{ for all } i\ll0,\\
\Longleftrightarrow\widehat{\Tor}_{i}^R(L,M^{\dagger})&=0  \text{ for all } i\ll0,
\end{array}\]
Set $S=R/\fc$ and consider the exact sequence (\ref{ttt}), 
$0\rightarrow\Hom_S(M,S)\rightarrow F\rightarrow N\rightarrow0,$ where $F$ is a free $S$--module. 
By \cite[Proposition 2.8]{CJ1}, the above exact sequence induces the following doubly infinite long exact sequence
\begin{equation}\tag{\ref{t44}.2}
\cdots\rightarrow\widehat{\Tor}_i^R(L,\Hom_S(M,S))\rightarrow\widehat{\Tor}_i^R(L,F)\rightarrow\widehat{\Tor}_i^R(L,N)\rightarrow\cdots.
\end{equation}
As $\pd_R(S)<\infty$, $\widehat{\Tor}_{i}^R(L,F)=0$ for all $i\in\ZZ$ by \cite[Lemma 2.7]{CJ1}. Note that $M^{\dagger}\cong\Hom_{S}(M,S)$ by \cite[Corollary]{G}. It follows from the exact sequece (\ref{t44}.2) that 
\[\begin{array}{rl}\tag{\ref{t44}.3}
\widehat{\Tor}_i^R(L, M^{\dagger})&=0 \text{ for all } i\ll0\\
\Longleftrightarrow\widehat{\Tor}_{i}^R(L,\Hom_S(M,S))&=0 \text{ for all } i\ll0\\
\Longleftrightarrow\widehat{\Tor}_{i}^R(L,N)&=0  \text{ for all } i\ll0,
\end{array}\]
Now the assertion is clear by (\ref{t44}.1) and (\ref{t44}.3).
\end{proof}
%%%%%%%%%%%%%%%%%%%%%%%%%%%%%%%%%%%%%%%%%%%%%%%%%%%%%%%%%%%%%%%%%%%%%%%%%%%%%%%%%%%%%
Similarly, one can prove the following result by using Proposition \ref{pr1} and Corollary \ref{c3}.
\begin{prop}\label{t5}
Let $R$ be a Gorenstein ring and let $\fa$ be a perfect Gorenstein ideal of $R$. Assume that $M\in\cm^n(R)$, $X\in\cm^m(R)$ and that $M\underset{\fa}{\sim}N$.
The following statements hold:
\begin{enumerate}[(i)]
\item{$\widehat{\Ext}^i_R(M,X)=0$ for $i\gg0$ if and only if $\widehat{\Ext}^i_R(X^{\dagger},N)=0$ for $i\gg0$.}
\item{$\widehat{\Ext}^i_R(M,X)=0$ for $i\ll0$ if and only if $\widehat{\Ext}^i_R(X^{\dagger},N)=0$ for $i\ll0$.}
\item{$\widehat{\Ext}^i_R(X,M)=0$ for $i\in\ZZ$ if and only if $\widehat{\Ext}^i_R(X^{\dagger},N)=0$ for $i\in\ZZ$.}
\end{enumerate}
\end{prop}
%%%%%%%%%%%%%%%%%%%%%%%%%%%%%%%%%%%%%%%%%%%%%%%%%%%%%%%%%%%%%%%%%%%%%%%%%%%%
Let $R$ be a Gorenstein ring and let $\fa$ be a perfect Gorenstein ideal of $R$. Assume that $M$ is a Cohen-Macaulay $R$--module and that $M\underset{\fa}{\sim}N$. In \cite[Corollary 1.3]{Pu}, Puthenpurakal proved that $\Ext^i_R(M,M)=0$ for all $i\gg0$ if and only if $\Ext^i_R(N,N)=0$ for all $i\gg0$. In the following, it is shown that $\Ext^i_R(M,M)\cong\Ext^i_R(N,N)$ for all $i>\G-dim_R(M)$.
\begin{thm}\label{t6}
Let $R$ be a Gorenstein ring and let $\fa$, $\fb$ be perfect Gorenstein ideals of $R$. Assume that $M\in\cm^n(R)$, $X\in\cm^m(R)$ and that $M\underset{\fa}{\sim}N$, $X\underset{\fb}{\sim}Y$. Then the following statements hold:
\begin{enumerate}[(i)]
\item{$\widehat{\Ext}^i_R(M,X)\cong\widehat{\Ext}^{i+m-n}_R(Y,N)$ for all $i\in\ZZ$.}
\item{$\widehat{\Ext}^i_R(M,M)\cong\widehat{\Ext}^i_R(N,N)$ for all $i\in\ZZ$.}
\item{If $M$ and $X$ are self-linked, then\\
$\widehat{\Ext}^i_R(M,X)\cong\widehat{\Ext}^{i+m-n}_R(X,M) \text{ for all } i\in\ZZ.$}
\end{enumerate}
\end{thm}
\begin{proof}
We only need to prove part (i), since parts (ii) and (iii) follow immediately from (i).
Set $X^{\ddagger}=\Ext^{m}_R(X,R)$ and $M^{\dagger}=\Ext^{n}_R(M,R)$. By Proposition \ref{pr1} and Lemma \ref{l5},
\[\begin{array}{rl}\tag{\ref{t6}.1}
\widehat{\Ext}^i_R(M,X)&\cong\widehat{\Tor}_{-i-1+n}^R(M^{\dagger},X)\\
&\cong\widehat{\Tor}_{-i-1+n}^R(X,M^{\dagger})\\
&\cong\widehat{\Ext}^{i+m-n}_R(X^{\ddagger},M^{\dagger}),
\end{array}\]
for all $i\in\ZZ$.
By (\ref{ttt}), we obtain the following exact sequences:
\begin{equation}\tag{\ref{t6}.2}
0\rightarrow\Hom_{R/\fa}(M,R/\fa)\rightarrow P\rightarrow N\rightarrow0,
\end{equation}
\begin{equation}\tag{\ref{t6}.3}
0\rightarrow\Hom_{R/\fb}(X,R/\fb)\rightarrow Q\rightarrow Y\rightarrow0,
\end{equation}
where $P$ is a free $R/\fa$--module and $Q$ is a free $R/\fb$--module. As $\fa$ and $\fb$ are perfect, $P$ and $Q$ have finite projective dimension and so $\widehat{\Ext}^i_R(Y,P)=0=\widehat{\Ext}^i_R(Q,M^{\dagger})$ for all $i\in\ZZ$ by Theorem \ref{tate}(i).
Note that $X^{\ddagger}\cong\Hom_{R/\fb}(X,R/\fb)$ and $M^{\dagger}\cong\Hom_{R/\fa}(M,R/\fa)$ by \cite[Corollary]{G}. It follows from the exact sequences (\ref{t6}.2), (\ref{t6}.3) and Theorem \ref{tate}(iv), (v) that
\[\begin{array}{rl}\tag{\ref{t6}.4}
\widehat{\Ext}^i_R(X^{\ddagger},M^{\dagger})&\cong\widehat{\Ext}^{i+1}_R(Y,M^{\dagger})\\
&\cong\widehat{\Ext}^{i}_R(Y,N),
\end{array}\]
for all $i\in\ZZ$. Now the assertion is clear by (\ref{t6}.1) and (\ref{t6}.4).
\end{proof}
%%%%%%%%%%%%%%%%%%%%%%%%%%%%%%%%%%%%%%%%%%%%%%%%%%%%%%%%%%%%%%%%%%%%%%%%%%%%%%%%%%%%%%%%%%%%%%
The following result is an immediate consequence of Theorem \ref{t6} and \cite[Theorem 5.9]{AM}.
\begin{cor}\cite[Proposition 10]{MS}
Let $R$ be a Gorenstein ring and let $\fa$ be a perfect Gorenstein ideal of $R$. Assume that $M$ and $N$ are $R$--modules and that $M\underset{\fa}{\sim}N$. Then $M$ is perfect if and only if $N$ is.
\end{cor}
%%%%%%%%%%%%%%%%%%%%%%%%%%%%%%%%%%%%%%%%%%%%%%%%%%%%%%%%%%%%%%%%%%%%%
Recall that an $R$--module $M$ is said to be \emph{horizontally self-linked} if $M\cong\lambda_RM$ (see \cite[Definition 7]{MS}).
Let $R$ be a complete Kleinian singularity. Equivalently, $R$ is the complete two-dimensional hypersurface ring $\mathbb{C}[|x, y, z|]/(f )$, where $f$ is one of Arnold's simple singularities. Recall that such singularities are indexed by the
Dynkin diagrams of types A,D, and E:
\[\begin{array}{rl}
A_n & f=x^2+y^2+z^{n+1},  n\geq1,\\
D_n & f=x^2+y^2z+z^{n-1},  n\geq4,\\
E_6 & f=x^2+y^3+z^4,\\
E_7 & f=x^2+y^3+yz^3,\\
E_8 & f=x^2+y^3+z^5.
\end{array}\]
In \cite[Theorem 3]{MS}, Martsinkovsky and Strooker proved that every indecomposable non-projective
maximal Cohen-Macaulay module over a complete Kleinian singularity is horizontally self-linked.
As a consequence of Theorem \ref{t6}, we obtain the following surprising result.
\begin{cor}
Let $R$ be a complete Kleinian singularity and let $M$, $N$ be $R$--modules.
Then $\widehat{\Ext}^i_R(M,N)\cong\widehat{\Ext}^{i}_R(N,M)$  for all  $i\in\ZZ$.
\end{cor}
\begin{proof}
Without lose of generality we may assume that $M$ and $N$ are indecomposable non-projective maximal Cohen-Macaulay $R$--modules by Remark \ref{rem} and \cite[Proposition 5.7]{AM}. Therefore, $M$ and $N$ are self-linked by \cite[Theorem 3]{MS}.
Now the assertion is clear by Theorem \ref{t6}.
\end{proof}
%%%%%%%%%%%%%%%%%%%%%%%%%%%%%%%%%%%%%%%%%%%%%%%%%%%%%%%%%%%%%%%%%%%%%%%%%%%%%%%%%%%%%%%%%%%%%%%
Recall that an ideal $I$ is said to be almost complete intersection if $\mu(I)\leq\hte(I)+1$, where $\mu(I)$ is the minimal number of generators of $I$.
\begin{cor}\label{c5}
Let $R$ be a Gorenstein local ring and let $M$ be an $R$--module. Assume that $\fa$ is a Cohen-Macaulay almost complete intersection ideal which is a generic complete intersection. Then $\Ext^i_R(M,R/\fa)=0$ for all $i\gg0$ if and only if $\Tor_i^R(M,R/\fa)=0$ for all $i\gg0$.
\end{cor}
\begin{proof}
By \cite[Proposition 1]{Sc}, there exists a perfect Gorenstein ideal $\fc$ such that $\fa$ and $\fb=(\fc:\fa)$
are linked and $\fb$ is a Gorenstein ideal. 
In other words, $R/\fa\cong\lambda_{R/\fc}(R/\fb)$ and $R/\fb\cong\lambda_{R/\fc}(R/\fa)$ by \cite[Proposition 1]{MS}.
It follows from Proposition \ref{t44}(i) that
$$\Ext^i_R(M,R/\fa)=0 \text{ for all } i\gg0$$
$$\Longleftrightarrow\Tor_{i}^R(M,R/\fb)=0  \text{ for all } i\gg0$$
By Proposition \ref{pr2} and (\ref{t}),
$$\Longleftrightarrow\widehat{\Ext}^{i}_R(M,R/\fb)=0 \text{ for all } i\gg0$$
By Proposition \ref{t44}(i),
$$\Longleftrightarrow\widehat{\Tor}_{i}^R(M,R/\fa)=0 \text{ for all } i\gg0.$$
\end{proof}
%%%%%%%%%%%%%%%%%%%%%%%%%%%%%%%%%%%%%%%%%%%%%%%%%%%%%%%%%%%
 The following is an immediate consequence of Corollary \ref{c5}, \cite[Corollary 4.4]{AY} and \cite[Proposition 3.2]{CJ}.
\begin{cor}
Let $R$ be a Gorenstein local ring and let $\fa$ be a Cohen-Macaulay almost complete intersection ideal which is a generic complete intersection.
Assume that either $\CI_R(R/\fa)<\infty$ or $R$ is AB ring. Then $\pd_R(R/\fa)<\infty$ if and only if $\Tor_i^R(R/\fa,R/\fa)=0$ for all $i\gg0$
\end{cor}
%%%%%%%%%%%%%%%%%%%%%%%%%%%%%%%%%%%%%%%%%%%%%%%%%%%%%%%%%%%%%%%%%%%%%%%%%%%%%%%%%%%%%%%%%%
In the following, it is shown that the stable Betti and Bass numbers are preserved under evenly linkage.
\begin{thm}\label{t2}
Let $R$ be a Gorenstein local ring of dimension $d$ and let $\fa$, $\fb$ be perfect Gorenstein ideals of $R$. Assume that $M\in\cm^n(R)$ and that $M\underset{\fa}{\sim}N$. Then the following statements hold:
\begin{enumerate}[(i)]
\item{$\widehat{\beta}_i^R(M)=\widehat{\mu}^{m+i-1}_R(N)$ for all $i\in\ZZ$, where $m=\depth_R(M)$.
In particular, $\beta_i^R(M)=\mu^{m+i-1}_R(N)$ for all $i>n+1$.}
\item{If $M$ is evenly linked to $L$, $M\underset{\fa}{\sim}N\underset{\fb}{\sim}L$, then $\widehat{\beta}_i^R(M)=\widehat{\beta}_i^R(L)$ and
$\widehat{\mu}^i_R(M)=\widehat{\mu}^i_R(L)$ for all $i\in\ZZ$. In particular,
$\beta_i^R(M)=\beta_i^R(L)$ for all $i>n$ and $\mu^i_R(M)=\mu^i_R(L)$ for all $i>d$.}
\end{enumerate}
\end{thm}
\begin{proof}
We only need to prove part (i), since part (ii) follow immediately from (i). Set $S=R/\fa$. Note that $M^{\dagger}\cong\Hom_{S}(M,S)$ by \cite[Corollary]{G}. Consider the following exact sequence $0\rightarrow\Hom_S(M,S)\rightarrow F\rightarrow N\rightarrow0,$ where $F$ is a free $S$--module (see \ref{ttt}).
It follows from the above exact sequence and \cite[Theorem 9.1]{AM} that $\widehat{\beta}_i^R(M)=\widehat{\beta}_{n-i}^R(N)$. Now the assertion is clear by Corollary \ref{c1} and Theorem \ref{G}(ii).
\end{proof}
We end the paper by the following result.
\begin{cor}
Let $R$ be a Gorenstein local ring and let $\fa$ be a Cohen-Macaulay almost complete intersection ideal which is a generic complete intersection.
Then $\widehat{\beta}_i^R(R/\fa)=\widehat{\mu}^{m+i-2}_R(R/\fa)$ for all $i\in\ZZ$, where $m=\dim(R/\fa)$.
In particular, $\beta_i^R(R/\fa)=\mu^{m+i-2}_R(R/\fa)$ for all $i>\G-dim_R(R/\fa)+2$.
\end{cor}
\begin{proof}
By \cite[Proposition 1]{Sc}, there exists a perfect Gorenstein ideal $\fc$ such that $\fa$ and $\fb=(\fc:\fa)$
are linked and $\fb$ is a Gorenstein ideal. By Theorem \ref{t2} and Corollary \ref{c6},
\[\begin{array}{rl}\tag{\ref{t6}.1}
\widehat{\beta}_i^R(R/\fa)&=\widehat{\mu}^{m+i-1}_R(R/\fb)\\
&=\widehat{\beta}_{i-1}^R(R/\fb)\\
&=\widehat{\mu}^{m+i-2}_R(R/\fa),
\end{array}\]
for all $i\in\ZZ$.
\end{proof}
%%%%%%%%%%%%%%%%%%%%%%%%%%%%%%%%%%%%%%%%%%%%%%%%%%%%%%%%%%%%%%%%%%%%%%%%%%%%%%%%%%%%%%%%%
%%-----------------------------------------------------------
%%% ----------------------------------------------------------------------
%%% ----------------------------------------------------------------------
\bibliographystyle{amsplain}
%%% ----------------------------------------------------------------------
%%% ----------------------------------------------------------------------
%%% ----------------------------------------------------------------------

\end{document}